\documentclass[12pt]{amsart} 
\usepackage[mathscr]{eucal}
\parskip=\smallskipamount

\newenvironment{matr}[1]{\left(\begin{array}{#1}}{\end{array} \right)}
\newcommand{\bmatr}[1]{\begin{matr}{#1}}
\newcommand{\ematr}{\end{matr}}

\newtheorem{theorem}{Theorem}[section]
\newtheorem{lemma}[theorem]{Lemma}

\newtheorem{proposition}[theorem]{Proposition}

\title{Research Problem: The Completion Number of a Graph}
\author{M. Bakonyi and T. Constantinescu}
\date{}

\address{Department of Mathematics\\
  Georgia State University\\
  Atlanta, GA 30303-3083}
\email{mbakonyi@cs.gsu.edu}

\address{Department of Mathematics \\
  University of Texas at Dallas \\
  Box 830688, Richardson, TX 75083-0688}
\email{tiberiu@utdallas.edu}

\begin{document}
\openup 2\jot
\advance\baselineskip by 1pt

\begin{abstract}
\noindent
Motivated by the remarkable interplay between chordal graphs and 
matrix algebra, we associate to each graph a so-called completion 
number that might encode some aspects of that interplay. We show that 
this number is not trivial, and we ask for a graph theoretic characterization
of those graphs with a given completion number. 
\end{abstract}

\maketitle

\bigskip
\section{Introduction}

A {\it partial matrix} is 
a matrix in which some entries are specified and the 
others are left as free complex variables. A {\it completion} of a partial
matrix is a conventional matrix obtained by the specification of all
free variables. We call a matrix $A=\{a_{ij}\}_{i,j=1}^n$ {\it
partial Hermitian}
if $a_{ii}$ is a specified real number for every $i=1,2,\cdots ,n$, $a_{ij}$
is specified if and only if $a_{ji}$ is specified, and $\bar{a}_{ji}=a_{ij}$.
In this paper we are only interested in Hermitian completions of Hermitian
partial matrices. A partial Hermitian matrix $A$ is called {\it partial
positive (semidefinite)} if all specified principal submatrices of $A$
are positive.

Let $G=(V,E)$ be an undirected graph in which $V$ is the finite
set of vertices and $E\subseteq \{(x,y)|x,y\in V, x\neq y
\}$ is the set of edges.
We shall assume throughout that $V=\{1,2,\cdots ,n\}$ for a positive
integer $n$. A subset $K\subseteq \{1,2,\cdots ,n\}$ is called a {\it clique} 
if $(x,y)\in E$ for every $x,y\in K$, $x\neq y$.
The graph $G=(\{1,2,\cdots ,n\},E)$ is said
to be {\it the graph of the Hermitian partial matrix} $A=\{a_{ij}\}_{i,j=1}^n$
if $a_{ij}$, $i\neq j$ is specified if and only if $(i,j)\in E$.
If $K$ is a clique of $G$, we denote by $A(K)$ the (specified) principal
submatrix $\{a_{ij}\}_{i,j\in K}$.

The {\it inertia} of an $n\times n$ Hermitian matrix $B$ is the triple 
$(i_+(B), i_-(B), i_0(B))$, in which $i_+(B)$ (resp. $i_-(B)$) is the 
number of positive (resp. negative) eigenvalues of $B$ counting
multiplicities, and $i_0(B)=n-i_+(B)-i_-(B)$ is the number of zero
eigenvalues of $B$. The following is an immediate consequence of the
so-called interlacing inequalities.

\begin{lemma}
\label{um}
Let $B$ be an $(n-1)\times (n-1)$ principal submatrix of an
$n\times n$ Hermitian matrix $M$. Then

\[
i_{\pm}(B)\le i_{\pm}(M)\le i_{\pm}(B)+1.
\]

\end{lemma}

Let $A$ be a partial Hermitian matrix with graph $G$ and let $K$ be a 
clique of $G$. By Lemma \ref{um}, for every Hermitian completion $M$
of $A$ we have that $i_{\pm}(M)\ge i_{\pm}(A(K))$. We define 

\[
i(A,-)=\max i_-(A(K))
\]

and

\[
i(A;0,-)=\max [i_0(A(K))+i_-(A(K))],
\]
in which the maximum is taken over all maximal cliques $K$ of $G$.

An undirected graph $G$ is called {\it chordal} if it has no minimal
simple circuits with $4$ or more edges. The following is the main result 
of \cite{JR}.

\begin{theorem}
\label{t1}
Let $G=(\{1,2,\cdots ,n\}, E)$ be a chordal graph, and let
$A$ be a partial Hermitian matrix with graph $G$. Then:

\noindent
(a) There exists a Hermitian completion $M$ of $A$ such that

\[
i_-(M)+i_0(M)=i(A;0,-).
\]

\noindent
(b) If, in addition, $A(K)$ is nonsingular for every maximal clique
$K$ of $A$, then there exists a nonsingular Hermitian completion $M$ for
which $i_-(M)=i(A,-)$. If $A$ is real, then in either statement (a)
or (b) M can be taken real as well.
\end{theorem}

The following well-known result in \cite{GJSW}
is a corollary of Theorem \ref{t1} for the case when $i(A,-)=0$ 

\begin{theorem}
\label{t2}
Every partial positive
matrix with a chordal graph admits a positive completion. 
\end{theorem}

In \cite{GJSW}
it was also shown that given any nonchordal graph $G$ there exists a partial
positive matrix $A$ with graph $G$ such that $A$ does not admit a positive
completion.

It appears as a potentially interesting question to ask 
what is the minimum number of negative eigenvalues one can always obtain
in case of a given nonchordal graph. 

\bigskip

\noindent
{\bf Definition.} Let $G$ be an undirected graph. The {\it
completion number} of $G$ is the smallest number $k$ such that every 
partial positive matrix $A$ with graph $G$ has a completion $M$ with
$i_-(M)\le k$.

\bigskip

The following is our proposed Research Problem.

\bigskip

\noindent
{\bf Question:} Given a graph $G$, how can one find its completion
number?

\bigskip
More specifically, it would be of interest to find graph theoretic
characterizations of those graphs with given completion number.
By the results in \cite{GJSW}, the completion number of a graph $G$ is $0$
if and only if $G$ is chordal. The completion number does not distinguish
the lenght of a chordless cycle with $n\ge 4$ vertices. More preciselly, 
if $G$ is the chordless $4$ cycle with vertices $\{1,2,3,4\}$, by adding 
the edge $13$ to $G$ we obtain the graph $G'$ which 
has the maximal cliques $\{1,2,3\}$ and $\{1,3,4\}$ that are not
cliques in $G$. If $A$ is a partial positive matrix with graph $G$
then we can specify the $13$ and $31$ entries to obtain a partial
Hermitian matrix $A'$ such that $A'(\{1,2,3\})$ is positive while
by Lemma \ref{um}, $i_-(A'(\{1,3,4\}))\le 1$.
Since $G'$ is chordal, by Theorem \ref{t1} $A'$ has a Hermitian 
completion $M$ with $i_-(M)\le 1$. This shows that the completion
number of $G$ is $1$. Similarly, the completion number of any 
chordless cycle with $n\ge 4$ vertices is $1$.

In fact, we can show that any positive integer is the completion number
of a certain graph. First, we show this for $n=2$.

\begin{proposition}
\label{p1}
Let $G_2=(V_2,E_2)$ be the graph with $V_2=\{1, 2, \cdots ,8\}$ and
$E_2$ defined by its complement as $E_2^c=\{13,24,57,68\}$.
Then the completion number of $G_2$ is $2$.
\end{proposition}

\begin{proof}
Let 

\[
A=\bmatr{cccccccc}
1 & 1 & ? & -1 & 0 & 0 & 0 & 0\\
1 & 1 & 1 & ? & 0 & 0 & 0 & 0\\
? & 1 & 1 & 1 & 0 & 0 & 0 & 0\\
-1 & ? & 1 & 1 & 0 & 0 & 0 & 0\\
0 & 0 & 0 & 0 & 1 & 1 & ? & -1\\
0 & 0 & 0 & 0 & 1 & 1 & 1 & ?\\
0 & 0 & 0 & 0 & ? & 1 & 1 & 1\\
0 & 0 & 0 & 0 & -1 & ? & 1 & 1\\
\end{matr} ,
\]
where ? denotes an unspecified entry. 
The graph of $A$ is $G_2$ and $A$ is partial positive. 
Since every Hermitian completion of both the upper-left and lower-right
$4\times 4$ principal submatrices of $A$ (these matrices were first
used as counterexamples in \cite{GJSW})
have at least one negative
eigenvalue, it follows that the completion number of $G_2$ is at
least $2$.

To prove that the completion number
is at most $2$, let $A$ be a partial
positive matrix with graph $G_2$. Consider $G_2^{(1)}$ be the
graph obtained by adding the edge $13$ to $G_2$. The following
is the list of the maximal cliques of $G_2^{(1)}$ that are not cliques in $G_2$:
$\{1,2,3,5,6\}$,
$\{1,2,3,6,7\}$,
$\{1,2,3,7,8\}$,
$\{1,2,3,5,8\}$,
$\{1,3,4,5,6\}$,
$\{1,3,4,6,7\}$,
$\{1,3,4,7,8\}$, and
$\{1,3,4,5,8\}$.
By Theorem \ref{t2}, we can specify the $13$ and the $31$ entires of $A$ to
obtain a partial Hermitian matrix $A_1$ with graph $G_2^{(1)}$ such that
$A_1(\{1,2,3,5,6\})$
is positive definite. By Lemma
\ref{um} the principal submatrices of
$A_1$ corresponding to the next $7$ maximal cliques have at most one negative
eigenvalue.

Let $G_2^{(2)}$ be the graph obtained by adding the
edge $24$ to $G_2^{(1)}$. The following is 
the list of the maximal cliques of $G_2^{(2)}$ that are not cliques in
$G_2^{(1)}$:
$\{1,2,3,4,5,6\}$,
$\{1,2,3,4,6,7\}$,
$\{1,2,3,4,7,8\}$, and
$\{1,2,3,4,5,8\}$.
By Theorem \ref{t1}
we can specify the $24$ and $42$ entries of $A_1$ to obtain 
a partial Hermitian matrix $A_2$ with graph $G_2^{(2)}$ such that
$i_-(A(\{1,2,3,4,
6, 7\}))\le 1$. By Lemma \ref{um},  
$i_-(A_2(\{1,2,3,4,5,6\}))\le 1$, 
while the principal submatrices of $A_2$ corresponding to the cliques
$\{1,2,3,4,7,8\}$ and
$\{1,2,3,4,5,8\}$ have at most two negative eigenvalues.

Let $G_2^{(3)}$ be the graph obtained by adding the
edge $57$ to $G_2^{(2)}$. Then 
$\{1,2,3,4,5,6,7\}$ and
$\{1,2,3,4,5,7,8\}$ the maximal cliques of $G_2^{(3)}$ that
are not cliques in $G_2^{(2)}$. By Theorem \ref{t1}, we can specify the
$57$ and $75$ entries of $A_2$ to obtain a partial Hermitian matrix
$A_3$ with graph $G_2^{(3)}$ such that 
$i_-(A_3(\{
1,2,3,4,5,7,8\}))\le 2$.
By Lemma \ref{um}, $i_-(A_3(
\{1,2,3,4,5,6,7\}))\le 2$.
Since $G_2^{(3)}$ is chordal, we can finally specify the
$68$ and $86$ entries of $A_3$ to obtain a Hermitian matrix
with at most two negative eigenvalues.
This completes the proof of the proposition.
\end{proof}

The general case is settled by a similar construction.

\begin{proposition}
\label{p2}
Let $G_n=(V_n,E_n)$ be the graph with $V_n=\{1,2,\cdots ,4n\}$
and $E_n$ defined by its complement as $E_n=\{(4k+1, 4k+3), (4k+2,
4k+4):
k=0,1, \cdots ,n-1\}$. Then the completion number of $G_n$ is $n$.
\end{proposition}

\begin{proof}
A straightforward generalization of the proof of Proposition \ref{p1}.
\end{proof}

There are many other questions that can be asked in connection with the 
completion number. For instance, 
is there any connection between the completion number of a graph
and its order as defined in \cite{AHMR}? Any connection to \cite{BJL}?

\end{document}